\documentclass[11pt]{article}

\usepackage{amscd}      
\usepackage{amssymb}
\usepackage[intlimits]{amsmath}
\usepackage{xypic}      
\LaTeXdiagrams          
\usepackage[all,v2]{xy}
\xyoption{2cell}
\UseAllTwocells
\xyoption{frame}
\CompileMatrices

\usepackage{theorem}

\newtheorem{prop}{Proposition}[section]
\newtheorem{lem}[prop]{Lemma}

\newtheorem{them}[prop]{Theorem}

\theorembodyfont{\upshape}

\newtheorem{defn}[prop]{Definition}

\newtheorem{rmk}[prop]{Remark}

\newenvironment{pf}{\begin{trivlist}\item[]{\sc Proof.}}%
            {\nolinebreak $\Box$ \end{trivlist}}
            {\nolinebreak $\Box$ \end{trivlist}}

\newcommand{\noprint}[1]{}

\renewcommand{\tilde}{\widetilde}

\newcommand{\toto}{\rightrightarrows}

\newcommand{\lcom}{_{\scriptscriptstyle\bullet}}

\newcommand{\zz}{{\mathbb Z}}

\newcommand{\nn}{{\mathbb N}}

\newcommand{\cc}{{\mathbb C}}
\newcommand{\rr}{{\mathbb R}}

\newcommand{\aA}{{\cal A}}

\newcommand{\dD}{{\cal D}}

\newcommand{\eE}{{\cal E}}

\newcommand{\fF}{{\cal F}}
\newcommand{\kK}{{\cal K}}
\newcommand{\lL}{{\cal L}}
\newcommand{\mM}{{\cal M}}
\newcommand{\hH}{{\cal H}}

\newcommand{\Cl}{\mathop{\rm Cliff}}

\newcommand{\ldiag}[1]%
       {\makebox[0cm]{${\scriptstyle#1}\downarrow\phantom{\scriptstyle#1}$}}
\newcommand{\ldiagup}[1]%
       {\makebox[0cm]{${\scriptstyle#1}\uparrow\phantom{\scriptstyle#1}$}}
\newcommand{\rdiag}[1]%
       {\makebox[0cm]{$\phantom{\scriptstyle#1}\downarrow{\scriptstyle#1}$}}
\newcommand{\sediagr}[1]%
       {\makebox[0cm]{$\phantom{\scriptstyle#1}\searrow{\scriptstyle#1}$}}
\newcommand{\nediagr}[1]%
       {\makebox[0cm]{$\phantom{\scriptstyle#1}\nearrow{\scriptstyle#1}$}}
\newcommand{\rdiagup}[1]%
       {\makebox[0cm]{$\phantom{\scriptstyle#1}\uparrow{\scriptstyle#1}$}}
\newcommand{\swdiag}[1]%
       {\makebox[0cm]{$\phantom{\scriptstyle#1}\swarrow{\scriptstyle#1}$}}
\newcommand{\sediag}[1]%
       {\makebox[0cm]{${\scriptstyle#1}\searrow\phantom{\scriptstyle#1}$}}
\newcommand{\nediag}[1]%
       {\makebox[0cm]{${\scriptstyle#1}\nearrow\phantom{\scriptstyle#1}$}}

\newcommand{\doublearrowstack}[2]%
 {{{{\scriptstyle#1}\atop{\textstyle\longrightarrow}}\atop{{\textstyle\longright
arrow}\atop{\scriptstyle#2}}}}
\newcommand{\rightleftarrowstack}[2]%
 {{{{\scriptstyle#1}\atop{\textstyle\longrightarrow}}\atop{{\textstyle\longlefta
rrow}\atop{\scriptstyle#2}}}}
\newcommand{\leftrightarrowstack}[2]%
 {{{{\scriptstyle#1}\atop{\textstyle\longleftarrow}}\atop{{\textstyle\longrighta
rrow}\atop{\scriptstyle#2}}}}

\newcommand{\overtoparrow}%
{\makebox[0cm]{\beginpicture
\setcoordinatesystem units <.8cm,.4cm> point at 0 0
\setplotarea x from -3 to 3, y from 0 to 1
\setquadratic
\plot -3 0 0 1 3 0 /
\put{\vector(3,-1){0}}[Bl] at 3 0
\endpicture}}

\newcommand{\underbottomarrow}%
{\makebox[0cm]{\beginpicture
\setcoordinatesystem units <.8cm,.4cm> point at 0 0
\setplotarea x from -3 to 3, y from 0 to 1
\setquadratic
\plot -3 1 0 0 3 1 /
\put{\vector(3,1){0}}[Bl] at 3 1
\endpicture}}

\newcommand{\ses}[5]%
{0\longrightarrow#1\stackrel{#2}{ \longrightarrow}#3\stackrel{#4}{
\longrightarrow}#5\longrightarrow0}

\newcommand{\dt}[6]%
{#1\stackrel{#2}{\longrightarrow}#3 \stackrel{#4}{\longrightarrow}#5
\stackrel{#6}{\longrightarrow} #1[1]}

\newcommand{\cat}[1]%
{(\mbox{\rm #1})}

\newcommand{\gm}{\Gamma}

\newcommand{\be }{\begin{eqnarray*}}
\newcommand{\ee }{\end{eqnarray*}}

\def\gpd{\,\lower1pt\hbox{$\longrightarrow$}\hskip-.24in\raise2pt
             \hbox{$\longrightarrow$}\,}

\newcommand\0{^{(0)}}
\newcommand\ho{\hat\otimes}
\newcommand{\clif}{{\mathbb{C}}\ell_1}
\newcommand{\zzz}{{\mathbb{Z}}/2{\mathbb{Z}}}
\newcommand{\gExt}{\mathop{\widehat{{\rm Ext}}}\nolimits}

\title{Twisted $K$-theory and Poincar\'e duality}

\begin{document}
\sloppy
\maketitle
\begin{abstract}
Using methods of $KK$-theory, we generalize Poincar\'e $K$-duality
to the framework of twisted $K$-theory.
\end{abstract}

\section*{Introduction}
In \cite{CS84}, Connes and Skandalis showed,
using Kasparov's $KK$-theory, that given
a compact manifold $M$, the $K$-theory of $M$ is isomorphic to the
$K$-homology of $TM$ and vice-versa. It is well-known to experts that
a similar result holds in twisted $K$-theory, although this is
apparently written nowhere in the literature. In this paper, using
Kasparov's more direct approach, we show that given any (graded)
locally trivial bundle $\aA$ of elementary $C^*$-algebras over $M$,
the $C^*$-algebras of continuous sections $C(M,\aA)$ and $C(M,\aA^{op}\ho
\Cl(TM\otimes \cc))$ are $K$-dual to each other. When $\aA=M$ is the
trivial bundle, we recover Poincar\'e duality between $C(M)$ and
$C_\tau(M):=C(M,\Cl(TM\otimes\cc))$ \cite{Kas88}, which is equivalent
to Poincar\'e duality between $C(M)$ and $C_0(TM)$ since $C_\tau(M)$
and $C_0(TM)$ are $KK$-equivalent to each other.

\section{Preliminaries}
In this paper, we will assume that the reader is familiar with the
language of groupoids (although this is not crucial in the proof
of the main theorem concerning Poincar\'e duality).

We just recall the definition of a generalized morphism
(see e.g. \cite{HS87}), since it is
used at several places. Suppose that $G\toto G^{(0)}$ and
$\gm\toto\gm^{(0)}$ are two Lie groupoids. Then a generalized morphism
from $G$ to $\gm$ is given by a space $P$, two maps
$G^{(0)}\stackrel{\tau}{\leftarrow} P\stackrel{\sigma}{\to} \gm^{(0)}$,
a left action of $G$ on $P$ with respect to $\tau$, a right action of $\gm$
on $P$ with respect to $\sigma$,
such that the two actions commute, and $P\to G^{(0)}$
and is a right $\gm$-principal bundle. The set of isomorphism classes
of generalized morphisms from $G$ to $\gm$ is denoted by $H^1(G,\gm)$.
There is a category whose objects are Lie groupoids and arrows are
isomorphism classes of generalized morphisms; isomorphisms in this category
are called Morita equivalences.

If $f:G\to\gm$ is a map such that $f(gh)=f(g)f(h)$ whenever $g$ and $h$
are composable, then $f$ is called a (strict) morphism. Then $f$
determines a generalized morphism $P_f=G^{(0)}\times_{\gm^{(0)}}\gm$.
Two strict morphisms $f$ and $f'$ determine the same element of
$H^1(G,\gm)$ if there exists $\lambda:G^{(0)}\to\gm$ such that
$f'(g)=\lambda(t(g)) f(g)\lambda(s(g))^{-1}$.

Finally, we recall that any element of $H^1(G,\gm)$ is given
by the composition of a Morita equivalence with a strict morphism.

\section{Graded twists and twisted $K$-theory}
In this section, we review the basic theory of twisted $K$-theory in
the graded setting, sometimes in more detail than some other references like
\cite{AS04,FHT,BS}. This is a probably well-known and straightforward
generalization of the ungraded case as developed e.g. in
\cite{AS04,TXL04}, hence we will omit most proofs.

\subsection{Graded Dixmier-Douady bundles}
Let $\mM\toto\mM\0$ be a Lie groupoid (more generally, most of the
theory below is still valid for locally compact groupoids having a Haar
system). The reader who is not interested in equivariant $K$-theory
may assume that $\mM=\mM\0=M$ is just a compact manifold.

A \emph{graded Dixmier-Douady bundle} of parity 0 (resp. of parity 1) $\aA$
over $\mM\toto \mM\0$ is a locally trivial bundle of $\zzz$-graded
$C^*$-algebras over $\mM\0$, endowed with a continuous action of
$\mM$, such that for all $x\in \mM\0$, the fiber $\aA_x$ is isomorphic
to the $\zzz$-graded algebra $\kK(\hat{H}_x)$ of compact operators over a
$\zzz$-graded Hilbert space $\hat{H}_x$ (resp. to the $\zzz$-graded
algebra $\kK(H_x)\oplus \kK(H_x)\cong \kK(H_x)\ho\clif$, where $H_x$ is
some Hilbert space and $\clif$ is the first complex Clifford algebra).
Beware that $H_x$ or $\hat{H}_x$ does not necessarily depend
continuously on $x$. Of course, the action of $\mM$ is required to
preserve the degree. The usual theory of graded twists \cite{AS04}
corresponds to even graded D-D bundles (i.e. D-D bundles of parity 0),
but our slightly more general definition allows to cover Clifford bundles as
well: if $E\to M$ is a Euclidean vector bundle of dimension $d$, then
$\Cl(E\otimes\cc)\to M$ is a graded D-D bundle of parity ($d\mbox{ mod }2$).

Denote by $\hat H$ the graded Hilbert space $H^0\oplus H^1$,
where $H^i=\ell^2(\nn)$, and $\hat H_\mM=L^2(\mM)\otimes \hat H$,
where $L^2(\mM)$ is the $\mM$-equivariant $\mM\0$-Hilbert module
obtained from $C_c(\mM)$ by completion with respect to the scalar product
$$\langle \xi,\eta\rangle (x)=\int_{g\in\mM^x}\overline{\xi(g)}\eta(g)\,
\lambda^x(dg).$$

Two graded D-D bundles $\aA$ and $\aA'$
are said to be Morita equivalent if (they have
the same parity and) $\aA\ho\kK(\hat H_\mM)\cong \aA'\ho\kK(\hat H_\mM)$.
The set of Morita equivalence classes of graded D-D bundles forms a
group $\widehat{Br}_*(\mM)=\widehat{Br}_0(\mM)\oplus \widehat{Br}_1(\mM)$,
the graded Brauer group of $\mM$. The sum of $\aA$ and $\aA'$ is
$\aA\ho\aA'$ (note that the parities do add up), and the opposite
$\aA^{op}$ of $\aA$ is the bundle whose fibre at $x\in \mM\0$ is the
conjugate algebra of $\aA_x$. In other words, $(\aA^{op})_x=
\kK(\hat{H}^*_x)$ (resp. $(\aA^{op})_x=\kK(H_x^*)\oplus \kK(H_x^*)$)
in the even (resp. odd) case.

Moreover, $\sigma_i:\aA\mapsto \aA\ho\clif$ is an isomorphism from
$\widehat{Br}_i(\mM)$ to $\widehat{Br}_{1-i}(\mM)$
such that $\sigma^2={\mathrm{Id}}$,
hence $\widehat{Br}_*(\mM)\cong
\widehat{Br}_0(\mM)\times \zzz$. Therefore, to study
$\widehat{Br}_*(\mM)$ it suffices to study $\widehat{Br}_0(\mM)$.
\par\medskip

Let us examine the relation between the graded Brauer group
$\widehat{Br}_0(\mM)$ and bundles of projective unitary operators.
Given any two graded Hilbert spaces $\hat{H}_1$ and $\hat{H}_2$, we denote
by $\hat{U}(\hat{H}_1,\hat{H}_2)$ the set of unitary operators from $\hat{H}_1$
to $\hat{H}_2$ which are homogeneous of degree 0 or 1, and
by $P\hat{U}(\hat{H}_1,\hat{H}_2)$ its quotient by $S^1$.
When $\hat{H}_1=\hat{H}_2$, these sets will be denoted by
$\hat{U}(\hat{H}_1)$ and $P\hat{U}(\hat{H}_1)$.

The set $H^1(\mM,P\hat{U}(\hat{H}))$ is actually an abelian monoid:
given two generalized morphisms $f_1$ and $f_2$ from
$\mM$ to $P\hat{U}(\hat{H})$, the composition of
$(f_1,f_2):\mM\to P\hat{U}(\hat{H})\times P\hat{U}(\hat{H})$
with the morphism $(u,v)\mapsto u\ho v$ is a generalized morphism
from $\mM$ to $P\hat{U}(\hat{H}\ho\hat{H})\cong
P\hat{U}(\hat{H})$.

Note that $H^1(\mM,\hat{U}(\hat{H}))$ is not a monoid, since
given two morphisms $f_1$, $f_2:\gm\to \hat{U}(\hat{H})$
(with $\gm$ Morita equivalent to $\mM$), the map
$f:g\mapsto f_1(g)\ho f_2(g)$ is not a morphism since

\begin{equation}\label{eqn:f1otimesf2}
f(gh)=(-1)^{|f_2(g)|\,|f_1(h)|}f(g)f(h)
\end{equation}

On the other hand, if we restrict to degree 0 operators,
i.e. if we consider $\hat{U}(\hat{H})^0 \cong U(H^0)\times U(H^1)$,
then $H^1(\mM,\hat{U}(\hat{H})^0)$ is again a monoid.

The sequence
$$H^1(\mM,\hat{U}(\hat{H})^0) \to
  H^1(\mM,P\hat{U}(\hat{H})) \to
  \widehat{Br}_0(\mM)\to 0,$$
where the first map is the quotient map and the second is
$P\mapsto P\times_{P\hat{U}(\hat{H})} \kK(\hat{H})$,
is canonically split-exact (the proof is analogue to \cite{TXL04}),
and the splitting identifies $\widehat{Br}_0(\mM)$ with
$H^1(\mM,P\hat{U}(\hat{H}))_{stable}=
\{[P]\vert\; [P]=[P+P_0]\}$, where
  $P_0=P\hat{U}(L^2(\hat{H}_\mM),\hat{H})$.
\par\medskip

Let us recall the relation with the ordinary Brauer group $Br(\mM)$
of Morita equivalence classes of ungraded D-D bundles. Recall that
$Br(\mM)\cong H^2(\mM\lcom,\underline{S^1})$, where $\underline{S^1}$ is the
sheaf of smooth $S^1$-valued functions, and that $Br(\mM)\cong
H^3(\mM\lcom,\zz)$
when $\mM\toto\mM\0$ is a proper groupoid, for instance, when
$(\mM\toto\mM\0)=(M\times G\toto M)$ is the crossed-product of a
manifold by a proper action of a Lie group $G$.

There is a split exact sequence \cite{AS04,FHT}
\begin{equation}\label{eqn:graded-brauer}
0\to Br(\mM)\to \widehat{Br}_0(\mM)\to
H^1(\mM,\zzz)\to 0.
\end{equation}

Indeed, from the exact sequence
$1\to P\hat{U}(\hat{H})^0 \to P\hat{U}(\hat{H})\to \zzz\to 0$,
we get an exact sequence
$0\to H^1(\mM,P\hat{U}(\hat{H})^0)_{stable}
 \to H^1(\mM, P\hat{U}(\hat{H}))_{stable}\to H^1(\mM,\zzz)\to 0$.
Moreover, there is an isomorphism
$H^1(\mM,P\hat{U}(\hat{H})^0)_{stable}
\cong H^2(\mM,\underline{S}^1)\cong Br(\mM)$ (this is analogue to the
fact that $H^1(\mM,PU(H))\cong Br(\mM)$).
\par\medskip
Furthermore, in the decomposition $\widehat{Br}_0(\mM)\cong
H^1(\mM,\zzz)\times Br(\mM)$, the sum becomes
$$(\delta_1,\aA_1)+(\delta_2,\aA_2)=
(\delta_1+\delta_1,\aA_1+\aA_2+\delta_1\cdot\delta_2),$$
where $\delta_1\cdot\delta_2$ is the element of $H^2(\mM,\underline{S}^1)$
corresponding to the cocycle $(\delta_1\cdot\delta_2)(g,h)
=(-1)^{\delta_2(g)\delta_1(h)}$. This can be seen by direct checking
using (\ref{eqn:f1otimesf2}) (see \cite[Proposition~2.3]{AS04} for a different
explanation).

\subsection{Graded $S^1$-central extensions}
\begin{defn}
A graded $S^1$-central extension of a groupoid $\gm\toto \gm^{(0)}$
is a central extension $S^1\to\tilde{\gm}\stackrel{\pi}{\to}\gm$,
together with a groupoid morphism $\delta:\gm\to\zzz$.
\end{defn}

One defines the sum of two graded central extensions $(\tilde{\gm}_1,\delta_1)$
and $(\tilde{\gm_2},\delta_2)$ as $(\tilde{\gm},\delta)$,
where $\delta(g)=\delta_1(g)+\delta_2(g)$ and
$\tilde{\gm}=(\tilde{\gm}_1\times_\gm \tilde{\gm}_2)/S^1
=\{(g_1,g_2)\in \tilde{\gm}_1\times \tilde{\gm}_2\vert\; \pi_1(g_1)
=\pi_2(g_2)\}/\sim$, where $\sim$ is the equivalence relation
$(g_1,g_2)\sim (g'_1,g'_2)\iff \exists \lambda\in S^1,\;
(g'_1,g'_2)=(\lambda g_1,\lambda^{-1}g_2)$.

The multiplication for the groupoid $\tilde{\gm}$ is
$(\tilde{g}_1,\tilde{g}_2)(\tilde{h}_1,\tilde{h}_2)=
(-1)^{\delta_2(g)\delta_1(h)}(\tilde{g}_1\tilde{h}_1,\tilde{g_2}\tilde{h}_2)$,
where $g=\pi_i(\tilde{g}_i)$, $h=\pi_i(\tilde{h}_i)$.

Note that the set of isomorphism classes of graded $S^1$-central extensions
of $\gm$ forms an abelian group. To see that the product is commutative,
if $\tilde{\gm}'=(\tilde{\gm}_2\times_\gm\tilde{\gm}_1)/S^1$ is
endowed with the product $(\tilde{g}_2,\tilde{g}_1)(\tilde{h}_2,\tilde{h}_1)
=(-1)^{\delta_1(g)\delta_2(h)}(\tilde{g}_2\tilde{h}_2,\tilde{g}_1\tilde{h}_1)$,
then
\begin{eqnarray*}
\tilde{\gm}&\to&\tilde{\gm'}\\
(\tilde{g}_1,\tilde{g}_2)&\mapsto& (-1)^{\delta_1(g)\delta_2(g)}
(\tilde{g}_2,\tilde{g}_1)
\end{eqnarray*}
is a $S^1$-equivariant isomorphism.
\par\medskip

To see that $(\tilde{\gm},\delta)$ has an inverse, let $\tilde{\gm}^{op}$
be equal to $\tilde{\gm}$ as a set, but the $S^1$-principal bundle
structure is replaced by the conjugate one, and the product $*_{op}$
in $\tilde{\gm}^{op}$ is
$$\tilde{g} *_{op} \tilde{h} = (-1)^{\delta(g)\delta(h)}
\tilde{g}\tilde{h}.$$
Then
\begin{eqnarray*}
\gm\times S^1 &\to& (\tilde{\gm}\times_\gm\tilde{\gm}^{op})/S^1\\
(g,\lambda)&\mapsto& [\lambda \tilde{g},\tilde{g}]
\end{eqnarray*}
is an isomorphism ($\tilde{g}\in\tilde{\gm}$ is any lift of
$g\in\gm$).
\par\medskip

Let us define the group $\gExt(\mM,S^1)$. Consider the collection
of triples $(S^1\to\tilde{\gm}\to\gm,
\delta,P)$ where $(S^1\to\tilde{\gm}\to \gm,\delta)$ is a graded central
extension and $P$ is a Morita equivalence from $\gm\to \mM$.
Two such triples
$(S^1\to\tilde{\gm}_1\to\gm_1,\delta_1,P_1)$ and $(S^1\to\tilde{\gm}_2\to\gm_2,
\delta_2,P_2)$ are said to be Morita equivalent if there exists a
Morita equivalence $\tilde{Q}:\tilde{\gm}_1\to\tilde{\gm}_2$
which is $S^1$-equivariant, such that the diagrams of isomorphism classes
of generalized morphisms
$$\xymatrix{
\gm_1\ar[r]^{[Q]}\ar[dr]_{[P_1]} & \gm_2\ar[d]^{[P_2]}\\
&\mM
}$$
and
$$\xymatrix{
\gm_1\ar[r]^{[Q]}\ar[dr]_{[\delta_1]} & \gm_2\ar[d]^{[\delta_2]}\\
&\zzz
}$$
commute, where $Q:\gm_1\to\gm_2$ is the Morita equivalence induced
by $\tilde{Q}$. Then the group $\gExt(\mM,S^1)$ is the quotient of the
collection of triples by Morita equivalence.
\par\medskip
Then $\widehat{Br}_0(\mM)\cong\gExt(\mM,S^1)$. Let us explain the map
$\widehat{Br}_0(\mM)\to\gExt(\mM,S^1)$.

Any element of $\widehat{Br}_0(\mM)$ is given by a bundle $\aA$ with fiber
$\kK(\hat{H})$, hence by a $P\hat{U}(\hat{H})$-principal bundle $P$ over
$\mM\toto \mM^{(0)}$, i.e. by a generalized morphism
$f:\mM\to P\hat{U}(\hat{H})$. Let $\eE$ be the graded central extension
$(S^1\to\hat{U}(\hat{H})\to P\hat{U}(\hat{H}),\delta)$, where
$\delta:P\hat{U}(\hat{H})\to \zzz$ is the degree map. Then $f^*\eE$ is
an element of $\gExt(\mM,S^1)$.

\begin{rmk}
Let $f_1^*\eE$, $f_2^*\eE\in\gExt(\mM,S^1)$. Then the sum of
$f_1^*\eE$ and $f_2^*\eE$ (where $f_i:\gm\to P\hat{U}(\hat{H})$ is
a strict morphism and $\gm\toto \gm^{(0)}$ is a groupoid which
is Morita-equivalent to $\mM$) is given by $f^*\eE$, where
\begin{eqnarray*}
f:\gm&\to& P\hat{U}(\hat{H}\ho\hat{H})\cong P\hat{U}(\hat{H})\\
g&\mapsto& f_1(g)\hat\otimes f_2(g).
\end{eqnarray*}
(Note that $f$ is indeed a homomorphism, since $f(gh)=
(-1)^{\delta(g)\delta(h)} f(g)f(h)$ agrees with $f(g)f(h)$
up to a scalar in $S^1$.)
\end{rmk}

\subsection{Twisted $K$-theory}
Let $\aA\to\mM^{(0)}$ be a graded D-D bundle over $\mM\toto \mM^{(0)}$.
We define $K_\aA^*(\mM)$ as $K_*(\aA\rtimes_{r}\mM)$,
the $K$-theory of the reduced crossed-product of the graded $C^*$-algebra
$\aA$ by the action of $\mM$. If $\aA$ corresponds to the graded central
extension $(S^1\to\tilde{\gm}\to\gm,\delta)$, then $K_\aA^*(\mM)$ is
isomorphic to $K_*(C^*_r(\tilde\gm)^{S^1})$,
where the $C^*$-algebra $C^*_r(\tilde\gm)^{S^1}$ is subalgebra of
$C^*_r(\tilde{\gm})$ which is the closure of
$$\{f\in C_c(\tilde{\gm})\vert\; f(\lambda g)=\lambda^{-1}f(g)\quad\forall
g\in\tilde{\gm}\;\forall\lambda\in S^1\}.$$
The $C^*$-algebra $C^*_r(\tilde{\gm})^{S^1}$ is considered as a $\zzz$-graded
$C^*$-algebra, using the grading automorphism
$$f\in C_c(\tilde{\gm})\mapsto (\gamma\mapsto f(\gamma)\delta(\gamma))
\in C_c(\tilde{\gm}).$$
\par\medskip

Note that it suffices to study $K^0_\aA(\mM)$, since
$K^1_\aA(\mM)=K^0_{\aA\ho\clif}(\mM)$.

\subsection{Example of manifolds}
Let $M$ be a manifold. Elements of $\gExt(M,S^1)$ are given by
an open cover $(U_i)_{i\in I}$, smooth maps $c_{ijk}:U_{ijk}=U_i\cap
U_j\cap U_k\to S^1$ and $\delta_{ij}:U_{ij}\to\zzz$ such that
$\delta_{ij}+\delta_{jk}=\delta_{ik}$ and
$c_{jkl}c_{ikl}^{-1}c_{ijl}c_{ijk}^{-1}=1$.

Let $\gm=\amalg_{i,j} U_{ij}$ and $\tilde{\gm}=\gm\times S^1$.
Define a product on $\tilde{\gm}$ by
$(x_{ij},\lambda)(x_{jk},\mu)=(x_{ik},\lambda\mu c_{ijk})$.
Then $\tilde{\gm}$ is a groupoid, and there is a central extension
$S^1\to\tilde{\gm}\to\gm$.

The sum of $(c^1,\delta^1)$ and $(c^2,\delta^2)$ is
$(c,\delta)$ where $\delta_{ij}=\delta^1_{ij}+\delta^2_{ij}$ and
$c_{ijk}=c^1_{ijk}c^2_{ijk}(-1)^{\delta_{ij}^2\delta^1_{jk}}$.
\par\medskip

Let us consider the particular case when $c=1$ is the trivial cocycle.
In that case, $C^*_r(\tilde{\gm})^{S^1}\cong C^*_r(\gm)$.
Let us compare this $\zzz$-graded $C^*$-algebra to the $\zzz$-graded
$C^*$-algebra $C_0(\tilde{M})$, where $\tilde{M}\to M$ is the
double cover determined by the cocycle $\delta$.
Let $P=(\amalg U_i)\times_M \tilde{M}$. Then $P$ is a $\zzz$-equivariant
Morita equivalence from $\gm$ to $\tilde{M}\rtimes \zzz$, hence
$K_\aA^*(M)=K_*(C^*_r(\tilde{\gm})^{S^1})=
K_*(C^*_r(\gm))=K_*(C_0(\tilde{M})\rtimes \zzz)=K_*(C_0(\tilde{M})
\ho\clif) = K_{*+1}(C_0(\tilde{M}))$ as in \cite[Remark~A.13]{FHT}.

\subsection{Twistings by Euclidean vector bundles}
Suppose that $E$ is a Euclidean vector bundle over $\mM\toto\mM^{(0)}$.
Then $E$ is given by an $O(n)$-principal bundle over $\mM\toto
\mM^{(0)}$, hence by a morphism $f:\gm\to O(n)$ together with a
Morita equivalence from $\gm$ to $\mM$.

Let $\eE$ be the graded $S^1$-central extension
$$S^1\to Pin^c(n)\to O(n),$$
where $Pin^c(n)=Pin(n)\times_{\{\pm 1\}} S^1$, and $\delta:O(n)\to\zzz$
is the map such that $\det A= (-1)^{\delta(A)}$. Then
$f^*\eE$ is a graded central extension of $\gm$, hence determines
an even graded D-D bundle $\aA_E$.

On the other hand, $\aA'_E={\mathrm{Cliff}}(E\otimes_\rr\cc)\to
\mM^{(0)}$ is another graded D-D bundle over $\mM\toto \mM^{(0)}$
which has the same parity as $\dim E$. We want to compare it to $\aA_E$.

We first need two lemmas.

\begin{lem}\label{lem:H3-trivial}
Let $G$ be a compact Lie group. Denote by $G_0$ its identity
component.
Assume that $G_0$ is simply connected and that $Br(G/G_0)=\{0\}$.
Then every central extension of $G$ by $S^1$ is split.
\end{lem}

\begin{pf}
Since $G$ is compact, every $S^1$-central extension is of finite order.
Let us recall the argument: given a central extension
$\eE=(S^1\to\tilde{G}\to G)$, let $V$ be a finite dimensional
representation of $\tilde{G}$ which is a
sub-representation of $\{f\in L^2(\tilde{G})\vert\;
f(\lambda g)=\lambda^{-1}f(g)\;\forall \lambda\in S^1,
\;\forall g\in \tilde{G}\}$. Let $d=\dim V$, $n=d!$ and $W=\Lambda^d V
\cong \cc$. Then the representation of $\tilde{G}$ in $W$
is a map $\tilde{G}\to U(W)\cong S^1$ which is a splitting of
$n\,\eE$, hence $\eE$ is of order at most $n$.
\par\medskip
Therefore, the extension $\eE$ comes from a central extension
$0\to \zz/n\zz\to\tilde{G}\to G\to 1$. Since $G$ is simply connected,
the central extension must be trivial as $\zz/n\zz$-principal bundle,
i.e. $\tilde{G}=G\times \zz/n\zz$, and the product on $\tilde{G}$
is given by $(g,\lambda)(h,\mu)=(gh,\lambda+\mu+ c(g,h))$ where $c:
G\times G\to \zz/n\zz$ is a 2-cocycle. Using connectedness of
$G_0$, $c$ must factor through $G/G_0\times G/G_0$, i.e.
the central extension is pulled back from a central extension
of $G/G_0$, which must be trivial by assumption.
\end{pf}

\begin{lem}\label{lem:ext-trivial}
Let $G$ be a Lie group and $G_0$ a normal subgroup containing the identity
component of $G$ such that $G_0$ has no nontrivial character and
that ${Br}(G_0)=\{0\}$.

If $\eE$ and $\eE'$ are $S^1$-central extensions whose restriction
to $G_0$ are isomorphic, then $\eE$ and $\eE'$ are isomorphic.
\end{lem}

\begin{pf}
After taking the difference of $\eE$ and $\eE'$, we may assume that
$\eE'$ is the trivial extension. Denote by $S^1\to\tilde{G}\to G$
the extension $\eE$. Let $g\mapsto \tilde{g}$ be
a splitting $G_0\to \tilde{G}$. Choose a family $(s_i)$ such that
$G=\amalg_i s_iG_0$, and for each $i$, choose a lift $\tilde{s}_i$
of $s_i$. Define then $\widetilde{s_ig}$ by $\tilde{s}_i\tilde{g}$.
By construction, $\tilde{\gamma h}=\tilde{\gamma}\tilde{h}$ for all
$(\gamma,h)\in G\times G_0$.

Next, define the 2-cocycle $c:G\times G\to S^1$ by
$\tilde{g}\tilde{h}=c(g,h)\tilde{gh}$. Let $c_{ij}=c(s_i,s_j)$.

For all $j$, let $\varphi_j:G_0\to S^1$ such that
$\tilde{s}_j^{-1}\tilde{g}\tilde{s}_j = \varphi_j(g)\widetilde{
s_j^{-1}gs_j}$. It is immediate to check that $\varphi_j$ is
a group morphism, hence $\varphi_j$ is trivial by assumption,
i.e. $\tilde{s}_j^{-1}\tilde{g}\tilde{s}_j = \widetilde{
s_j^{-1}gs_j}$. Multiplying on the right by $\tilde{h}$
and on the left by $\tilde{s}_i\tilde{s}_j$, we get
$\tilde{s_i}\tilde{g}\tilde{s}_j\tilde{h}=\tilde{s}_i\tilde{s}_j\widetilde{
s_j^{-1}gs_jh}=c_{ij}\widetilde{s_is_j}\widetilde{
s_j^{-1}gs_jh}=c_{ij}\widetilde{s_igs_jh}$, hence
$\widetilde{s_ig}\widetilde{s_jh}=c_{ij}\widetilde{s_igs_jh}$.
It follows that $c(s_ig,s_jh)=c_{ij}$, i.e. that $c$ factors
through $G/G_0\times G/G_0$. Since $Br(G/G_0)$ is trivial by assumption,
$c$ must be a coboundary. We conclude that $\eE$ is a split extension.
\end{pf}

\begin{prop}
$K_{*+\dim E,\aA_E}(\mM)=K_{*,\aA'_E}(\mM)=K_*(C_0(\mM^{(0)},
\mathop{\rm Cliff}(E\otimes_\rr\cc))\rtimes_r\mM)=K_*(C_0(E)\rtimes_r \mM)$.
\end{prop}

\begin{pf}
The last equality follows from the fact that $C_0(E)$ and $C_0(\mM^{(0)},
\mathop{\rm Cliff}(E\otimes_\rr\cc))$ are $\mM$-equivariantly $KK$-equivalent
\cite{Kas80}.
The second equality is just the definition.

To prove the first equality, let us suppose for instance that $n=\dim E$
is even, the proof for $n$ odd being analogous.
We have to compare the graded D-D bundle $\aA'_E$ with the graded
central extension $S^1\to\tilde{\gm}\to\gm$ which is pulled back from
$S^1\to Pin^c(n)\to O(n)$. By naturality, we can just assume that
$\gm=\mM=O(n)$, and that $E=\rr^n$ is endowed with the canonical action of
$O(n)$.

Then, $\mathop{\rm Cliff}(E\otimes_\rr\cc)=\cc\ell_n=\lL(\hat{\mathcal{H}})$,
where $\hat{\mathcal{H}}$ is the graded Hilbert space
$\cc^{2^{n/2-1}}\oplus \cc^{2^{n/2-1}}$.

Denote by $\alpha:O(n)\to P\hat{U}(\hat{\mathcal{H}})$ the canonical
action of $O(n)$ on $\cc\ell_n=\Cl(E_\cc)$. To show that the
central extension associated to the graded D-D bundle $\Cl(E_\cc)
\to\cdot$ is $(S^1\to Pin^c(n)\to O(n))$, it suffices to prove that there
exists a lifting $\tilde{\alpha}$:
$$\xymatrix{
S^1\ar[r]\ar[d] & S^1\ar[d]\\
Pin^c(n)\ar[d]\ar[r]^{\tilde{\alpha}} & \hat{U}(\hat{\hH})\ar[d]\\
O(n)\ar[r]^\alpha & P\hat{U}(\hat{\hH}).
}$$
For $n=2$ it is elementary to check that both vertical lines are
split extensions.
For $n\ge 3$, since $Spin(n)$ is simply connected and
compact, it has no nontrivial $S^1$-central extension
(see Lemma~\ref{lem:H3-trivial}), hence the pull-back of
$\eE=(S^1\to \hat{U}(\hat{\hH})\to P\hat{U}(\hat{H}))$
by the map $Spin(n)\to O(n)$ has a lift $\beta:Spin(n)\to \hat{U}
(\hat{\hH})$.

If $\beta(-1)=Id$ then $\beta$ induces a map $\bar\beta:SO(n)\to\hat{U}
(\hat{\hH})$, which means that the extension $S^1\to Spin^c(n)\to
SO(n)$ is split. It follows that $S^1\to Spin^c(3)\to SO(3)$ is
split, i.e. that there exists a morphism $\varphi:Spin^c(3)
=SU(2)\times_{\{\pm 1\}} S^1\to S^1$ such that $\varphi(\lambda g)=
\lambda \varphi(g)$ for all $\lambda\in S^1$ and all $g\in
SU(2)\times_{\{\pm 1\}} S^1$. Putting $\chi(g)=\varphi(g,1)$,
the morphism $\chi: SU(2)\to S^1$ satisfies $\chi(-1)=-1$.
Using simplicity of $SU(2)/{\{\pm 1\}}=SO(3)$,
it follows that $\chi$ is injective, which is impossible.
\par\medskip

It follows that $\beta(-1)=-Id$, hence $\beta$ induces a lift $\beta:
Spin^c(n)\to \hat{U}(\hat{\hH})$ which is $S^1$-equivariant.
This means that the restriction of $\eE$ to $SO(n)$ is
isomorphic to $S^1\to Spin^c(n)\to SO(n)$. To conclude that
the restriction of $\eE$ to $O(n)$ is isomorphic to
$S^1\to Pin^c(n)\to O(n)$, we apply Lemma~\ref{lem:ext-trivial}
to $G=O(n)$ and $G_0=SO(n)$.
\end{pf}

\section{Poincar\'e duality}
\subsection{Kasparov's constructions}
Let $M$ be a compact manifold (actually, Poincar\'e duality can be
generalized to arbitrary manifolds \cite{Kas88}, but in this paper we
confine ourselves to compact ones for simplicity).
We suppose that $M$ is endowed with a Riemannian metric which is
invariant by the action of a locally compact group $G$. Given any
vector bundle $\aA$ over any manifold $M$, we denote by $C_\aA(M)$
the space of continuous sections vanishing at infinity. We will
also write $C_\aA$ whenever there is no ambiguity.
We denote by $\tau$ the complexified cotangent bundle of $M$.
\par\medskip

In \cite{Kas88}, Kasparov constructed two elements
$$\theta\in KK_{M\rtimes G}(C(M),C(M)\otimes C_\tau(M))=
RKK_G(M;\cc,C_\tau(M))$$
and $D\in KK_G(C_\tau(M),\cc)$ (in this paper, we will use Le Gall's
\cite{Leg}
notation $KK_{M\rtimes G}(\cdot,\cdot)$ for
equivariant $KK$-theory with respect to the groupoid $M\rtimes G$,
rather than Kasparov's $RKK_G(M;\cdot,\cdot)$, but of course
both are equivalent).

Let us recall the construction of $\theta$ and $D$.

Let $H=L^2(\Lambda^* M)$, and
\begin{eqnarray*}
\varphi:C_\tau(M)&\to& \lL(H)\\
\omega&\mapsto& e(\omega)+e(\omega)^*,
\end{eqnarray*}
where $e(\omega)$ is the exterior multiplication,
and let $F=\dD(1+\dD)^{-1/2}$ where $\dD=d+d^*$.
Then $D=[(H,\varphi,F)]$.

Let us explain the construction of $\theta$. Denoting by $\rho$ the distance
function on $M$,
let $r>0$ be so small that for all $(x,y)$ in $U=\{(x,y)\in M\times M\vert\;
\rho(x,y)<r\}$, there exists a unique
geodesic from $x$ to $y$.

For every $C(M\times M)$-algebra $A$, we denote by $A_U$ the
$C^*$-algebra $C_0(U)A$. Then the element $\theta$ is defined
as $[(C_M\otimes C_\tau(M))_U,\Theta]$ where
$\Theta=(\Theta_x)_{x\in M}$,
$\Theta_x(y)=\frac{\rho(x,y)}{r}(d_y\rho)(x,y)\in T^*_yM
\subset \mathop{\rm Cliff}_\cc(T^*_yM)$.

\subsection{Constructions in twisted $K$-theory}
In this subsection, we construct an element $\theta^\aA\in
KK_{M\rtimes G}(C(M),C_\aA\ho C_{\aA\ho\aA^{op}})$ for any
graded D-D bundle $\aA$ over $(M\times G\toto M)$, i.e. for any
$G$-equivariant graded D-D bundle over $M$. We may assume that
$\aA$ is stabilized, i.e. that $\aA\cong \aA\ho \kK(\hat{H}\otimes L^2(G))$.
First, let us denote by $p_t(x,y)$ the geodesic segment joining $x$ to $y$
at constant speed ($0\le t\le 1$).

Using $p_t$, we see that $p_t:U\to M$ is a $G$-equivariant homotopy
equivalence. Unfortunately, this does not imply that $Br(U\rtimes G)$
and $Br(M\rtimes G)$ are isomorphic
for arbitrary $G$, hence we will make the following
\par\medskip
{\bf{Assumption.}}
In the sequel of this paper, and unless stated otherwise,
$G$ will be a \emph{compact} Lie group
acting smoothly on a compact manifold $M$.
\par\medskip
In that case, $H^2(U\rtimes G,\underline{S}^1)
\cong H^3(U\rtimes G,\zz) = H^3(\frac{U\times EG}{G},\zz)
\cong H^3(\frac{M\times EG}{G},\zz)\cong
H^2(M\rtimes G,\underline{S}^1)$.
As a consequence, there is a continuous, $G$-equivariant family
of isomorphisms
$$u_{t,x,y}:\aA_x\stackrel{\sim}{\to} \aA_{p_t(x,y)}.$$
Of course, the $u_t$'s are not unique,
but this will not be important as far as $K$-theory is concerned
as we will see.
\par\medskip
Consider the canonical Morita-equivalence $\hH_x$ between $\cc$ and
$\aA_x\ho\aA^{op}_x$. Let $\hH=(\hH_x)_{x\in M}$ be the corresponding
Morita equivalence between $C(M)$ and $C_{\aA\ho\aA^{op}}(M)$.
Then $[((C(M)\otimes C_\tau(M))_U\ho_{C(M)}\hH,\Theta\ho 1)]$
is an element of $KK_{M\rtimes G}(C(M),C_{\aA\ho\aA^{op}}(M)\ho C_\tau(M))$.

Now, using the map $u_{1,x,y}:\aA_x\stackrel{\sim}{\to}\aA_y$, we get a
Morita equivalence $\eE_{x,y}$
from $\aA_x^{op}$ to $\aA_y^{op}$, thus a Morita equivalence $\eE$
from $(C_{\aA^{op}}(M)\ho C(M))_U$ to $(C(M)\ho C_{\aA^{op}}(M))_U$.
Tensoring over $C(M\times M)$ with $C_\aA(M)\ho C_\tau(M)$, we get a Morita
equivalence $\eE'=(\eE'_{x,y})_{(x,y)\in U}$ from
$(C_{\aA\ho\aA^{op}}(M)\ho C_\tau(M))_U$ to
$(C_\aA(M)\ho C_{\tau\ho \aA^{op}}(M))_U$.

We then define $\theta^\aA$ as
\begin{eqnarray*}
\theta^\aA&=&[((C(M)\otimes C_\tau(M))_U\ho_{C(M)} \hH\ho_{C_{\aA\ho\aA^{op}}
\ho C_\tau}\eE',\Theta\ho 1)]\\
&&\in KK_{M\rtimes G}(C(M),C_\aA(M)\ho C_{\tau\ho \aA^{op}}(M)).
\end{eqnarray*}

\subsection{Twisted $K$-homology}
Given a $C^*$-algebra $A$ endowed with an action of a locally compact group
$G$, the $G$-equivariant $K$-homology of $A$, $K^*_G(A)$, is defined
by $KK^*_G(A,\cc)$. If $\aA$ is a $G$-equivariant graded D-D bundle over $M$,
we define $K^{G,\aA}_*(M)$ by $K^*_G(C_\aA(M))$.

\subsection{Maps between twisted $K$-theory and twisted $K$-homology}
Let $A$ and $B$ be two separable graded
$G$-$C^*$-algebras (recall that $G$ is
assumed to be a compact Lie group). We define two maps
\begin{eqnarray*}
\mu:KK_G(A,C_\aA(M)\ho B)&\to&KK_G(C_{\tau\ho \aA^{op}}\ho A,B)\\
\nu:KK_G(C_{\tau\ho \aA^{op}}\ho A,B)&\to& KK_G(A,C_\aA(M)\ho B).
\end{eqnarray*}

First, let us introduce some notations. Suppose that $M$ is a locally
compact space endowed with an action of a locally compact group.
If $A$, $B$ and $D$ are $G$-equivariant graded $C(M)$-algebras, and $\eE$ is a
$A$-$B$-$C^*$-bimodule, then Kasparov defined a
$A\ho_{C_0(M)} D$-$B\ho_{C_0(M)} D$-$C^*$-bimodule
$\sigma_{M,D}(\eE)$, and thus a ``suspension'' map $\sigma_{M,D}:
KK_{M\rtimes G}(A,B)\to KK_{M\rtimes G}(A\ho_{C_0(M)} D,B\ho_{C_0(M)} D)$.
There is also a suspension map
$\sigma_D:KK_{M\rtimes G}(A,B)\to KK_{M\rtimes G}(A\ho D,B\ho D)$
defined in a similar way.
\par\medskip

Given a $(A_1,B_1\ho D)$-
$C^*$-bimodule $\eE_1$ and a $(D\ho A_2,B_2)$-$C^*$-bimodule
$\eE_2$, the $(A_1\ho A_2,B_1\ho B_2)$-$C^*$-bimodule
$\eE_1\ho_D\eE_2$ is defined by
$\sigma_{A_2}(\eE_1)\ho_{A_2\ho B_1}\sigma_{B_1}(\eE_2)$.

We introduce a similar notation when all tensor products are replaced
by tensor products over a space $M$: given a $(A_1,B_1\ho_{C(M)} D)$-
$C^*$-bimodule $\eE_1$ and every $(D\ho_{C(M)} A_2,B_2)$-$C^*$-bimodule
$\eE_2$, $\eE_1\ho_D\eE_2$ is a
$(A_1\ho_{C(M)} A_2,B_1\ho_{C(M)} B_2)$-$C^*$-bimodule.
\par\medskip
Let us now define $\mu$ and $\nu$.
First, let us note that $KK_G(A,C_\aA(M)\ho B)$ is isomorphic to
$KK_{M\rtimes G}(C(M)\otimes A,C_\aA(M)\ho B)$.

The map $\mu$ is defined as the composition
\begin{eqnarray*}
KK_{M\rtimes G}(C(M)\otimes A,C_\aA(M)\ho B)
&\stackrel{\sigma_{M,C_{\tau\ho \aA^{op}}}}{\to}&
KK_{M\rtimes G}(C_{\tau\ho\aA^{op}}\ho \aA,C_{\tau\ho\aA^{op}\ho\aA}
\ho B)\\
&\stackrel{\cdot\ho \hH^{op}}{\to}&
KK_{M\rtimes G}(C_{\tau\ho \aA^{op}}\ho A,C_\tau\ho B)\\
&\stackrel{\otimes D}{\to}&
KK_G(C_{\tau\ho\aA^{op}}\ho A,B).
\end{eqnarray*}

The map $\nu$ is just $\theta^\aA\otimes\cdot:
KK_G(C_{\tau\ho\aA^{op}}\ho A,B)\to KK_{M\rtimes G}(C(M)\otimes A,
C_\aA\ho B)$.

\subsection{The main theorem}
\begin{them}\label{thm:main}
Let $G$ be a compact Lie group acting on a compact manifold $M$,
and let $A$ and $B$ be two graded separable $G$-$C^*$-algebras.
Let $\aA$ be a $G$-equivariant graded D-D bundle over $M$.

Then the maps $\mu$ and $\nu$ defined above are inverse to each other:
$$KK^*_G(A,C_\aA(M)\ho B)\cong KK^*_G(C_{\tau\ho\aA^{op}}(M)\ho A,B).$$
Replacing $\aA$ by $\tau\ho \aA^{op}$, we get:
$$KK_G^*(C_\aA(M)\ho A,B)\cong KK_G^*(A,C_{\tau\ho \aA^{op}}(M)\ho B).$$
In particular, for $A=B=\cc$ we get
\begin{eqnarray*}
K^*_{G,\aA}(M)&\cong& K_*^{G,\tau\ho\aA^{op}}(M)\\
K^{G,\aA}_*(M)&\cong& K^*_{G,\tau\ho \aA^{op}}(M).
\end{eqnarray*}
\end{them}

\begin{rmk}
This result (in the case when $G$ is the trivial group) is observed
in \cite[Section~7]{BMRS}.
\end{rmk}

\begin{rmk}
The map $\mu$ does not depend of the choice of the isomorphisms
$u_{t,x,y}$, hence $\nu$ doesn't either.
\end{rmk}

The rest of the paper is devoted to the proof of Theorem~\ref{thm:main}.

\subsection{Proof of $\mu\circ\nu=\mbox{Id}$}
For all $\alpha\in KK_G(C_{\tau\ho\aA^{op}}\ho A,B)$, we have
$$\mu\circ\nu(\alpha)=\sigma_{M,C_{\tau\ho\aA^{op}}}(\theta^\aA)
\ho_{C_{\aA^{op}\ho\aA}}\hH^{op}\otimes_{C_\tau}D
\otimes_{C_{\tau\ho\aA^{op}}\ho A} \alpha.$$
Thus, we need to prove that
$$\sigma_{M,C_{\tau\ho\aA^{op}}}(\theta^\aA)
\ho_{C_{\aA^{op}\ho\aA}}\hH^{op}\otimes_{C_\tau}D=1
\in KK_G(C_{\tau\ho \aA^{op}},C_{\tau\ho \aA^{op}}).$$
Consider the element $\sigma_{M,C_{\tau\ho\aA^{op}}}(\theta)\in
KK_{M\rtimes G}(C_{\tau\ho \aA^{op}},C_{\tau\ho\aA^{op}}\ho C_\tau)$.
Denote by
\begin{eqnarray*}
s:C_{\tau\ho\aA^{op}}\ho C_\tau &\stackrel{\cong}{\to} &
C_\tau\ho C_{\tau\ho \aA^{op}}\\
x\otimes y &\mapsto& (-1)^{\deg x.\deg y} y\otimes x
\end{eqnarray*}
the flip. Suppose proven that
\begin{equation}
\label{eqn:a}
\sigma_{M,C_{\tau\ho\aA^{op}}}(\theta)\otimes [s]=
\sigma_{M,C_{\tau\ho\aA^{op}}}(\theta^\aA)\ho_{C_{\aA^{op}\ho\aA}}
\hH^{op}.
\end{equation}
Then
\begin{eqnarray*}
\sigma_{M,C_{\tau\ho\aA^{op}}}(\theta^{\aA})\otimes_{C_{\aA^{op}\ho\aA}}
\hH^{op}\otimes_{C_\tau} D &=&
\sigma_{M,C_{\tau\ho \aA^{op}}}(\theta)\otimes_{C_\tau}D\\
&=& \sigma_{M,C_{\tau\hat\otimes \aA^{op}}}(\theta\otimes_{C_\tau}D)
=\sigma_{M,C_{\tau\ho\aA^{op}}}(1)=1
\end{eqnarray*}
Since $\theta\otimes_{C_\tau}D=1$ (from \cite[Theorem~4.8]{Kas88}).

We postpone the proof of (\ref{eqn:a}) until subsection~\ref{subsec:proofa}.

\subsection{Proof of $\nu\circ\mu=\mbox{Id}$}
For all $\alpha\in KK_{M\rtimes G}(C(M)\otimes A,C_\aA(M)\ho B)$,
we have
$$\nu\circ\mu(\alpha)=\theta^\aA\otimes_{C_{\tau\ho\aA^{op}}}
(\sigma_{M,C_{\tau\ho\aA^{op}}}(\alpha)\otimes_{C_{\aA^{op}\otimes \aA}}
\hH^{op}\otimes_{C_\tau}D).$$

Suppose shown that
\begin{eqnarray}\label{eqn:b}
\theta^\aA\otimes_{C_{\tau\ho\aA^{op}}}\sigma_{M,
C_{\tau\ho\aA^{op}}}(\alpha)\otimes_{C_{\aA^{op}\ho\aA}}\hH^{op}
&=&\alpha\otimes_{C_\aA}\sigma_{M,C_\aA}(\theta)\\
\nonumber
&&\in KK_{M\rtimes G}(C(M)\otimes A,C_\aA\ho C_\tau\ho B).
\end{eqnarray}
Then $\nu\circ\mu(\alpha)=\alpha\otimes_{C_\aA}
(\sigma_{M,C_\aA}(\theta)\otimes_{C_\tau}D)
=\alpha\otimes_{C_\aA}\sigma_{M,C_\aA}(\theta\otimes_{C_\tau}D)
=\alpha\otimes_{C_\aA}1=\alpha$.

We postpone the proof of (b) until subsection~\ref{subsec:proofb}

\subsection{Proof of (\ref{eqn:a})}
\label{subsec:proofa}
Recall the proof when $\aA$ is the trivial bundle \cite[Lemma~4.6]{Kas88}.
We want to show that $\sigma_{M,C_\tau}(\theta)$ is flip-invariant.
Denote by $p_t^*:T^*_{p_t(x,y)}M\hookrightarrow T^*_{(x,y)}U$ the
pull-back map induced by $p_t$, and let $q_t^*$ be the isometry
$q_t^*=p_t^*(p_tp_t^*)^{-1/2}$. We denote again by $q_t^*:
\Omega^1(M)\hookrightarrow \Omega^1(U)$ the corresponding map.
Then $q_t^*$ induces a map
$\varphi_t:C_\tau(M)\to \lL(C_\tau(U))$.

Let $\beta(s)=(p_{t(1-s)}(x,y),p_{t(1-s)+s}(x,y))\in U$, and
$$\Theta_t(x,y)=
\frac{\rho(x,y)}{r}\left\|\frac{d\beta}{ds}|_{s=1}\right\|^{-1}
\frac{d\beta}{ds}|_{s=1}\in T_{x,y}U.$$
Then $(C_\tau(U),\varphi_t,\Theta_t)_{0\le t\le 1}$
is a homotopy between $\sigma_{M,C_\tau}(\theta)$ and
$\sigma_{M,C_\tau}(\theta)\otimes [s]$.
\par\bigskip

Now, consider the general case.
$\sigma_{M,C_{\tau\ho\aA^{op}}}(\theta)$ is the Kasparov triple
$$((C_{\tau\ho\aA^{op}}\ho C_\tau)_U,\varphi,\Theta)$$
where $\varphi:C_{\tau\ho\aA^{op}}\to \lL((C_{\tau\ho\aA^{op}}\ho
C_\tau)_U)$ is the obvious map. Thus,
$$\sigma_{M,C_{\tau\ho\aA^{op}}}(\theta)\otimes [s]=
[((C_\tau\ho C_{\tau\ho\aA^{op}})_U,\varphi,\Theta_1)]$$
with $\theta_1=\frac{\rho(x,y)}{r} (d_x\rho)(x,y)\in
T_x^*M\subset \mathop{\rm Cliff}_\cc(T^*_xM)$,
while
$$\sigma_{M,C_{\tau\ho \aA^{op}}}(\theta^\aA)\otimes_{C_{\aA^{op}\ho
    \aA}}
\hH^{op}=[((C_{\tau\ho\aA^{op}}(M)\ho C_\tau(M))_U
\ho_{C_{\tau\ho\aA^{op}}\ho C_\tau}\eE'',\varphi',\Theta\otimes 1)]$$
where $\eE''$ is the Morita equivalence between $(C_{\tau\ho\aA^{op}}
\ho C_\tau)_U$ and $(C_\tau\ho C_{\tau\ho\aA^{op}})_U$
obtained from the Morita equivalence $\eE$ between
$p_0^* C_{\aA^{op}}=(C_{\aA^{op}}\otimes C(M))_U$ and
$p_1^* C_{\aA^{op}}=(C(M)\otimes C_{\aA^{op}})_U$.
\par\medskip
Let $\eE_t=(\eE_{x,y,t})_{(x,y)\in U}$
be the Morita equivalence between
$p_t^*C_{\aA^{op}}$ and $p_1^*C_{\aA^{op}}$ constructed in the same
    way as $\eE=\eE_0$.

Then
\begin{eqnarray*}
\sigma_{M,C_{\tau\ho\aA^{op}}}(\theta)\otimes [s]&=&
(C_\tau(U)\ho_{C_0(U)}\eE_1,\varphi,\Theta_1)\\
\sigma_{M,C_{\tau\ho\aA^{op}}}(\theta^\aA)\ho_{C_{\aA^{op}\ho\aA}}
\hH^{op}&=&(C_\tau(U)\ho_{C_0(U)}\eE_0,\varphi',\Theta).
\end{eqnarray*}

Let $\Theta_t$ as above. We produce a homotopy
$$(C_\tau(U)\ho_{C_0(U)}\eE_t,\psi_t,\Theta_t)$$
between those two elements.
Only $\psi_t:C_{\tau\ho\aA^{op}}\to\lL(C_\tau(U)\ho_{C_0(U)}\eE_t)$
remains to be defined. We need two compatible maps
\begin{eqnarray*}
\psi'_t:C_\tau&\to& \lL(C_\tau(U)\ho_{C_0(U)}\eE_t)\\
\mbox{and }\;\psi''_t:C_{\aA^{op}}&\to& \lL(C_\tau(U)\ho_{C_0(U)}\eE_t).
\end{eqnarray*}
The map $\psi'_t$ is just $\varphi_t\otimes 1$.
The map $\psi''_t$ is given by the composition
$$C_{\aA^{op}}\stackrel{p_t^*}{\to}
C_b(U,p_t^*\aA^{op})\to
\lL(\eE_t).$$

\subsection{Proof of (\ref{eqn:b})}
\label{subsec:proofb}

Let us first recall the proof when $\aA$ is trivial \cite[Lemma~
  4.5]{Kas88}.
We want to show that for all $\alpha\in KK_{M\rtimes G}(C(M)\otimes
A,C(M)\otimes B)$ we have
$$\alpha\otimes_{C(M)}\theta
=\theta\otimes_{C_\tau(M)}\sigma_{M,C_\tau(M)}(\alpha)
\in KK_{M\rtimes G}(C(M)\otimes A, C(M)\otimes C_\tau(M)\ho B).$$
Write $\alpha=[(E,T)]$ where $\overline{C(M,A)E}=E$ and $T$ is
$G$-continuous. Then both products can be written as
$$[(E\ho_{C(M)}(C(M)\otimes C_\tau(M))_U, \varphi_i,F_i)]$$
where $F_i$ is of the form $M_2^{1/2}(T\ho 1)+ M_1^{1/2}(1\ho\Theta)$
($i=0,1$), and where the map $C(M)\to
\lL((C(M)\otimes C_\tau(M))_U)$ used to define $\varphi_i$ is
$p_i^*$.

Since $p_0$ and $p_1$ are homotopic, $\varphi_0$ and $\varphi_1$ are
homotopic. One then constructs a homotopy between $F_0$ and $F_1$
using Kasparov's technical theorem as in \cite[Lemma~4.5]{Kas88}.
\par\bigskip
Let us now consider a general $G$-equivariant graded D-D bundle $\aA$
over $M$. Let $\alpha=[(E,T)]\in
KK_{M\rtimes G}(C(M)\otimes A,C_\aA\ho B)$ where
$\overline{C(M,A)\,E}=E$ and $T$ $G$-continuous.

We want to show that $\alpha\otimes_{C_\aA}\sigma_{M,C_\aA}(\theta)
=\theta^\aA\ho_{C_\tau\ho\aA^{op}}(\alpha)\ho_{C_{\aA^{op}\ho\aA}}
\hH^{op}$. Let us just explain the homotopy between the two modules,
the homotopy between the $F_i$'s being obtained using Kasparov's
technical theorem in the same way as in \cite[Lemma~4.5]{Kas88}.

The left-hand side is
\begin{equation}
\label{eqn:LHS}
E\ho_{C_\aA}(C_\aA\ho C_\tau)_U,
\end{equation}

and the right-hand side is
\begin{equation}
\label{eqn:RHS}
E\ho_{C_\aA}\sigma_{M,C_\aA}(\fF_1)\ho_{C_0(U)}
\hH'_1\ho_{C_{p_1^*(\aA\ho\aA^{op})}(U)} p_1^*\hH^{op},
\end{equation}
where we recall that $p_1:U\to M$ is the second projection $(x,y)\mapsto
y$. The $C(M)$-$(C(M)\otimes C_\tau)_U$-bimodule
$\fF_1$ is $(C(M)\otimes C_\tau)_U$, with the
left action of $C(M)$ on $\fF_1$ obtained via
$C(M)\stackrel{p_1^*}{\to} C_b(U)\to\lL((C(M)\otimes C_\tau)_U)$.

$\hH'_1$ is the Morita equivalence between $C_0(U)$ and $p_0^*C_{\aA}
\ho_{C_0(U)}
p_1^*C_{\aA^{op}}$ obtained by composing the Morita equivalence
$p_0^*\hH$ between $C_0(U)$ and $p_0^*(C_{\aA\ho\aA^{op}})$ with the
isomorphism $p_0^*\aA^{op}\cong p_1^*\aA^{op}$.

Using the map $p_t:U\to M$ instead of $p_1$, consider (with
obvious notations) the homotopy $E\ho_{C_\aA}
\sigma_{M,C_\aA}(\fF_t)\ho_{C_0(U)}\hH'_t\ho_{C_{p_t^*(\aA\ho\aA^{op})}(U)}
p_t^*\hH^{op}$.

For $t=1$, we get (\ref{eqn:RHS}).

For $t=0$, we get $E\ho_{C_\aA} (C_\aA\ho C_\tau)_U \ho_{C_0(U)}
p_0^*\hH\ho_{C_{p_0^*(\aA\ho\aA^{op})}(U)} p_0^*\hH^{op}$
where the right $C_{p_0^*(\aA\ho\aA^{op})}(U)$-structure on
$E\ho_{C_\aA} (C_\aA\ho C_\tau)_U \ho_{C_0(U)}
p_0^*\hH$ is defined as follows: $C_{p_0^*\aA}$ acts on $(C_\aA\ho
C_\tau)_U$ by the obvious action,
and $C_{p_0^*\aA^{op}}$ acts on $p_0^*\hH$. In other words,
it is the tensor product of (\ref{eqn:LHS}) with $\beta_\aA$ over
$C_\aA$, where $\beta_\aA$ is the
$C_\aA$-$C_\aA$-bimodule

$$\beta_\aA=(C_\aA\ho_{C(M)}\hH)\ho_{C_{\aA\ho\aA^{op}}}\hH^{op}.$$
In the expression above,
the right $C_{\aA\ho\aA^{op}}$-module structure on $C_\aA\ho_{C(M)}\hH$
is defined as follows: $\forall a\in C_\aA$, $\forall b\in C_{\aA^{op}}$
$\forall \xi\otimes\eta\in C_\aA\ho_{C(M)}\hH$,

$$(\xi\otimes\eta)\cdot (a\otimes b)=(-1)^{|\eta|\,|a|}\xi a \otimes \eta b.$$
To finish the proof, it remains to show that $\beta_\aA\cong C_\aA$.
Suppose for instance that $\aA$ is an even graded D-D bundle.
Let $x\in M$. Denoting by $\hat{H}_x$ a Hilbert space such that $\aA_x\cong
\kK(\hat{H}_x)$, we have $(\beta_\aA)_x =
[\kK(\hat{H}_x)\ho (\hat{H}_x^*\ho \hat{H}_x)]
\ho_{\aA_x\ho\aA_x^{op}} (\hat{H}_x\ho \hat{H}_x^*)$,
where in $\hat{H}_x\ho \hat{H}_x^*$, $\hat{H}_x$ (resp. $\hat{H}_x^*$) is considered as
a $\aA_x$-$\cc$ (resp. a $\aA_x^{op}$-$\cc$-bimodule, and in
$\hat{H}_x^*\ho \hat{H}_x$, $\hat{H}_x$ (resp. $\hat{H}_x^*$) is considered as
a $\cc$-$\aA_x^{op}$ (resp. a $\cc$-$\aA_x$-bimodule,
and the right $\aA_x\ho\aA_x^{op}$-module structure on $\kK(\hat{H}_x)\ho
(\hat{H}_x^*\ho \hat{H}_x)$ is $(\xi\otimes (\eta\otimes \zeta))\cdot
(a\otimes b)=(-1)^{|a|\,(|\xi|+|\eta|)} \xi a\otimes\eta\otimes\zeta b$.
It follows that $(\beta_\aA)_x\cong \hat{H}_x\ho \hat{H}_x^*$ is the natural
$\kK(\hat{H}_x)$-bimodule $\kK(\hat{H}_x)$.

\end{document}